\documentclass[11pt,twoside]{amsart}
\usepackage[utf8]{inputenc}
\usepackage{graphicx}
\usepackage{color}
\usepackage{amssymb}
\usepackage{amsmath}
\usepackage{amsthm}
\usepackage{amsfonts}
\usepackage{amsmath, amsfonts, amssymb}
\usepackage{tikz}
\usepackage{diagbox} 
\usepackage{array}

\newtheorem{theorem}{Theorem}[section]

\newtheorem{de}{Definition}

\newtheorem{theor}{Theorem}

\newtheorem{corollary}[theorem]{Corollary}

\numberwithin{equation}{section}

\begin{document}
	\title{On three-dimensional associative algebras}
	\author{U. Bekbaev$^1$,  I. Rakhimov$^2$}
	
	\thanks{{\scriptsize
			emails: $^1$uralbekbaev@gmail.com; isamiddin@uitm.edu.my$^2$ (corresponding author)}}
	\maketitle
	\begin{center}
		\address{$^{1}$Turin Polytechnic University in Tashkent, Tashkent, Uzbekistan,\\ $^2$Universiti Teknologi MARA (UiTM), Shah Alam, Malaysia and\\ V.I. Romanovsky Institute of Mathematics, Uzbekistan Academy of Sciences, Uzbekistan}
	\end{center}
	
	\begin{abstract}

This paper addresses the classification problem of associative algebras over arbitrary base fields. We present a list of three-dimensional associative algebras with canonical representatives of the isomorphism classes for fields of characteristic different from two and three. We compare our lists with the most recent classifications over the complex numbers and with the nilpotent case over arbitrary base fields in dimension three, adding some comments. Finally, we also provide a complete list of three-dimensional permutative algebras over a field of characteristic not two and three.

	\end{abstract}

	\section{Introduction}
The problem of classifying finite-dimensional algebras is one of the most challenging problems in algebra to date. For associative algebras in particular, the first attempts at classification were made at the end of the 19th century by B. Pierce.

The theory of finite-dimensional associative algebras is one of the oldest areas of modern algebra. It originates in the works of W.R. Hamilton, who discovered the famous quaternions, and A. Cayley, who developed the theory of matrices. Later, the structural theory of finite-dimensional associative algebras was studied by a number of mathematicians, notably B. Pierce, C.S. Pierce, W.K. Clifford, K. Weierstrass, R. Dedekind, C. Jordan, and F.G. Frobenius.

A new era in the development of the theory of finite-dimensional associative algebras began with the works of J. Wedderburn, who obtained several fundamental results: the description of the structure of semisimple algebras over a field, the theorem on lifting of quotients by the radical, the theorem on the commutativity of finite division rings, and others. Further progress was made by the German school of algebraists, headed by E. Noether, E. Artin, and R. Brauer. As a result, the theory of semisimple algebras reached its modern form, and most of the results were extended to certain classes of rings.

Further development of the theory of associative algebras took place in the 1980s, when many open problems that had remained unsolved since the 1930s were finally resolved. Nevertheless, researchers have continually returned to the classification problem of associative algebras, again and again, because of the wide range of applications and the deep connections of this class with other classes of algebras. A complete classification, or even a universal method of classification in any fixed dimension (even over the field of complex numbers), is still far from being achieved.

Many attempts have been made to obtain a complete classification of low-dimensional unital associative algebras, and a number of lists have been published. To mention just a few: the classification in dimensions up to four was obtained in 1935 by G. Scorza \cite{S} and later reproved by P. Gabriel using quivers in 1974 \cite{Gabriel}. In \cite{FP, FPP}, A. Fialowski, M. Penkava, and M. Phillipson provided a list of three-dimensional complex associative algebras by studying their deformations. G. Mazzola \cite{M} produced a classification in dimension five in studying the geometric classification, while M. Goze and A. Makhlouf \cite{GM} obtained a list in dimension six. B.Poonen classified  nilpotent commutative associative algebras of dimension $\leq 5$ over
algebraically closed fields \cite{P}. Note that all these results were obtained over the field of complex numbers.

Recently, we came across the paper \cite{KSTT}, where the authors claim to provide a complete classification of three-dimensional associative algebras over
$\mathbb{R}$ and $\mathbb{C}$. Having in hand the result of \cite{B}, which gives a complete classification of two-dimensional algebras over an arbitrary base field, we decided to apply it to the classification of three-dimensional associative algebras over any base field using a method called the extension method.

The essence of the method is as follows: we fix a two-dimensional subalgebra of the three-dimensional algebra we aim to construct and write down its matrix of structure constants. Taking these structure constants as unknowns and imposing the associativity condition yields a system of equations. For each fixed subalgebra, the resulting system of equations can be successfully solved, producing (in general) a redundant list of algebras containing that subalgebra. Since the classification of two-dimensional associative algebras over arbitrary fields was obtained in \cite{R}, we can use that as input. To make the list of three-dimensional algebras irredundant, we apply the automorphism groups of two-dimensional associative algebras over arbitrary fields. Both the construction of three-dimensional algebras and the elimination of redundancies turned out to be quite computationally intensive, and therefore we used Maple software to carry out these computations.

The organization of the paper is as follows.

In Section 2, we introduce some of the basic definitions and results which are used in the paper.

In Section 3, we present the list of all three-dimensional associative algebras over a field of characteristic is not two and three, followed by the description of the procedure that the list is complete and nonredundant.

In Section 4, we give the list of three-dimensional complex associative algebras obtained in \cite{KSTT}, compare the list with our list and provide comments. Here by Tables 1-3 we give the reasons why the extra algebras appeared in our list are not isomorphic those obtained in \cite{KSTT}.

In Section 5, we establish the correspondence between our list and W.A. De Graaf's list (see \cite{GR}) for three-dimensional nilpotent associative algebras over a field of characteristic not two and three.

Section 6 deals with the classification problem of three-dimensional permutative algebras. Since a permutative algebra is associative in \cite{Ivan} the authors gave their own classification of three-dimensional complex associative algebras. Therefore, first we compare all three lists and revise the list of three-dimensional complex associative algebras to make it consistent. Then we give a complete and nonredundant list of thee-dimensional (left and right) permutative algebras over a field of characteristic not two and three, and comparison with the list of left permutative algebras given in \cite{Ivan} over the field of complex numbers.

\section{Preliminaries}
All algebras in the paper are supposed to be finite-dimensional over a field $\mathbb{F}$.
\begin{de} \label{ISO} Two algebras $\mathbb{A}=(\mathbb{V},\cdot_{\mathbb{A}})$ and $\mathbb{B}=(\mathbb{V},\cdot_{\mathbb{B}})$ are called isomorphic
if there is an invertible linear map  $f:\mathbb{V}\rightarrow \mathbb{V} $
such that \[f(\mathrm{x}\cdot_{\mathbb{A}}
\mathrm{y})=f(\mathrm{x})\cdot_{\mathbb{B}} f(\mathrm{y}),\]
where $\mathrm{x}, \mathrm{y}\in \mathbb{A}.$\end{de}

\begin{de} An invertible linear map  $f:\mathbb{V}\rightarrow
\mathbb{V} $ is said to be an automorphism if \[f(\mathrm{x}\cdot
\mathrm{y})=f(\mathrm{x})\cdot f(\mathrm{y}),\] whenever
$\mathrm{x}, \mathrm{y}\in \mathbb{A}=(\mathbb{V},\cdot).$ \end{de}

The set of all automorphisms of an algebra $\mathbb{A}$ forms a group with respect to the composition operation and it is denoted by $Aut(\mathbb{A}).$

Let $\mathbb{A}=(\mathbb{V},\cdot)$ be an $n$-dimensional algebra over $\mathbb{F}$ and
$\mathbf{e}=(\mathrm{e}_1,\mathrm{e}_2,...,\mathrm{e}_n)$ be its basis. Then the bilinear map $\cdot:\mathbb{V} \times \mathbb{V}\longrightarrow \mathbb{V}$ is represented by a $n \times n^2$ matrix (called the matrix of structure constant, shortly MSC) $$A=\left(\begin{array}{ccccccccccccc}a_{11}^1&a_{12}^1&...&a_{1n}^1&a_{21}^1&a_{22}^1&...&a_{2n}^1&...&a_{n1}^1&a_{n2}^1&...&a_{nn}^1\\ a_{11}^2&a_{12}^2&...&a_{1n}^2&a_{21}^2&a_{22}^2&...&a_{2n}^2&...&a_{n1}^2&a_{n2}^2&...&a_{nn}^2 \\
...&...&...&...&...&...&...&...&...&...&...&...&...\\ a_{11}^n&a_{12}^n&...&a_{1n}^n&a_{21}^n&a_{22}^n&...&a_{2n}^n&...&a_{n1}^n&a_{n2}^n&...&a_{nn}^n\end{array}\right)$$ as follows:
$$\mathrm{e}_i \cdot \mathrm{e}_j=\sum\limits_{k=1}^n a_{ij}^k\mathrm{e}_k, \ \mbox{where}\ i,j=1,2,...,n.$$
This means that if a basis $\mathbf{e}$ is fixed then the algebra can be represented as a $(n \times n^2)$ matrix $A \in M_{n\times n^2}(\mathbb{F})$. Further once a basis is fixed we
do not make a difference between the algebra
$\mathbb{A}$ and its MSC $A$,  and work with matrices $A$.

The product on $\mathbb{A}$ with respect to the basis $\mathbf{e}$ is written as below
\begin{equation} \label{Product}
	\mathrm{x}\cdot \mathrm{y}=\mathbf{e} A(x\otimes y)
\end{equation}
 for any $\mathrm{x}=\mathbf{e}x,\mathrm{y}=\mathbf{e}y,$
	where $x=(x_1, x_2,...,x_n)^T,$ and  $y=(y_1, y_2,...,y_n)^T$ are column coordinate vectors of $\mathrm{x}$ and $\mathrm{y},$ respectively, $x\otimes y$ is the tensor (Kronecker)
 product of $x$ and $y$. Now and onward for the product ``$\mathrm{x}\cdot \mathrm{y}$'' on $\mathbb{A}$ we use the juxtaposition ``$\mathrm{x} \mathrm{y}$''.
\begin{de}  An algebra $\mathbb{A}$ is said to be associative if
for all $\mathrm{x}, \mathrm{y}, \mathrm{z} \in \mathbb{A}$ the following axiom is satisfied
\begin{equation} \label{AA}
\begin{array}{rll}
(\mathrm x  \mathrm y) \mathrm z &=& \mathrm x  (\mathrm y  \mathrm z).
\end{array}
\end{equation}
 \end{de}
Note the condition (\ref{AA}) in terms of MSC $A$ can be written as follows
$$\mathbf{e}A(A(x\otimes y)\otimes z))=\mathbf{e}A(x\otimes A(y\otimes z)),$$
i.e., an algebra $\mathbb{A}$ with MSC $A$ is associative if and only if
\begin{equation} \label{AA2}
\begin{array}{ccc}
A(A\otimes I)&=&A(I\otimes A),
\end{array}
\end{equation}
Note that in terms of MSC Definition \ref{ISO} can be rewritten as
\begin{equation}\label{ISO1}gB=Ag^{\otimes 2} \Longleftrightarrow B=g^{-1}Ag^{\otimes 2}\end{equation} and it produces a linear representation of $GL_n(\mathbb{F})$ on $M_{n\times n^2}(\mathbb{F})$ as follows
\begin{equation}\label{ISO1}\rho:(GL_n(\mathbb{F}),M_{n\times n^2}(\mathbb{F})) \longrightarrow M_{n\times n^2}(\mathbb{F}): \ \ \rho(g,A)=gA(g^{-1})^{\otimes 2}.\end{equation}
For $A \in M_{n\times n^2}(\mathbb{F})$ we introduce the vectors
$$Tr_1(A)=\left(\sum\limits_{j=1}^n a_{j1}^j,\ \sum\limits_{j=1}^n a_{j2}^j, ..., \sum\limits_{j=1}^n a_{jn}^j\right), \ \ Tr_2(A)=\left(\sum\limits_{j=1}^n a_{1j}^j,\ \sum\limits_{j=1}^n a_{2j}^j, ..., \sum\limits_{j=1}^n a_{nj}^j\right).$$
It is not hard to see, applying the contraction operation to the tensor $(a^i_{jk})_{i,j,k=1,...,n}$, that
 \begin{equation}Tr_1(gA(g^{-1})^{\otimes 2})=Tr_1(A)g^{-1},\ \ Tr_2(gA(g^{-1})^{\otimes 2})=Tr_2(A)g^{-1}.\end{equation}
 In particular, if $\lambda Tr_1(A)=Tr_2(A)$ is true then $\lambda Tr_1(gA(g^{-1})^{\otimes 2})=Tr_2(gA(g^{-1})^{\otimes 2})$ is valid as well.
  Therefore, we split $M(n)=M_{n\times n^2}(\mathbb{F})$ into the following disjoint subsets:
   	\begin{enumerate}
  		\item $M_2(n)=\left\{A \in M_{n\times n^2}(\mathbb{F}):\ \  \{Tr_1(A), \ Tr_2(A)\} \ \mbox{is linearly independent} \right\};$
  		\item $M_{1,\lambda}(n)=\left\{A \in M_{n\times n^2}(\mathbb{F}):\ \  \lambda Tr_1(A)= Tr_2(A) \ \mbox{and non of}\  Tr_1(A), Tr_2(A) \ \mbox{is zero}\right\};$
  		\item $M_{1,0}(n)=\left\{A \in M_{n\times n^2}(\mathbb{F}): \  \ Tr_1(A)\ \mbox{is nonzero and}\ Tr_2(A)\ \mbox{is zero }\right\};$
  		\item $M_{0,\infty}(n)=\left\{A \in M_{n\times n^2}(\mathbb{F}): \   Tr_1(A)  \ \mbox{is zero and}\ Tr_2(A)\ \mbox{is nonzero}\right\};$
  		\item $M_{0}(n)=\left\{ A \in M_{n\times n^2}(\mathbb{F}): \ \mbox{the both}\ Tr_1(A), Tr_2(A) \ \mbox{are zero} \right\}.$
  	\end{enumerate}
Let now turn back to three-dimensional case and
 $$A=
		 \begin{pmatrix}a_1& a_2& a_3& a_4&a_5 & a_6& a_7& a_8& a_9\\b_1& b_2& b_3& b_4&b_5 & b_6& b_7& b_8& b_9\\c_1& c_2& c_3& c_4&c_5 & c_6& c_7& c_8& c_9\end{pmatrix}$$
be MSC of $\mathbb{A}$. Then the trace vectors $Tr_1$ and $Tr_2$ are read as follows:
$$Tr_1=(a_1+b_4+c_7,a_2+b_5+c_8,a_3+b_6+c_9),\ \ Tr_2=(a_1+b_2+c_3,a_4+b_5+c_6,a_7+b_8+c_9).$$
The result that we are going to use from \cite{R} was given as follows.
\begin{theor} \label{Ass} Every non-trivial $2$-dimensional associative algebra over a field $\mathbb{F}$ $(Char(\mathbb{F})\neq 2)$
is
isomorphic to one of the following  presented by their matrices of structure constants pairwise non-isomorphic algebras:
\begin{enumerate}
\item  $As_{2}^1:=\left(\begin{array}{ccccc} 0&0&0&0\\
1&0&0&0
\end{array}\right),$
\item  $As_2^2:=$ $\left(\begin{array}{ccccc} 1&0&0&0\\
0&0&0&0
\end{array}\right),$
\item $As_2^3:=$
    $\left(\begin{array}{ccccc} 1&0&0&0\\
0&1&0&0
\end{array}\right),$
\item $As_2^4:=$ $\left(\begin{array}{ccccc} 1&0&0&0\\
0&0&1&0
\end{array}\right),$
\item $As_2^5(\alpha_4):=$ $\left(\begin{array}{ccccc} 1&0&0&\alpha_4\\
0&1&1&0
\end{array}\right)\simeq\begin{pmatrix}
 			1&0&0&a^2\alpha_4 \\
 			0& 1& 1& 0
 			\end{pmatrix},$ where $ \alpha_4 \in\mathbb{F},$ $a\in\mathbb{F}$ and $ a\neq 0.$
\end{enumerate}
\end{theor}
\section{List of three-dimensional associative algebras over a base field $\mathbb{F}$ $(Char(\mathbb{F})\neq 2,3)$}
Here are representatives of the isomorphism classes of nontrivial three-dimensional associative algebras over a base field $\mathbb{F}$ with the characteristic not two and three.
\begin{theor} \label{3-dim} Every three-dimensional associative algebra over a field $\mathbb{F}$ of characteristic not two and three is isomorphic to one the following pairwise non-isomorphic such algebras given by their MSC as follows
\begin{itemize}
\item Associative algebras with linearly independent traces $Tr_1$ and $Tr_2:$
	\begin{enumerate}
		\item
$As_2^2(3)(0)=\begin{pmatrix}1& 0& 0& 0& 0& 0& 0& 0& 0\\0& 1& 0& 1& 0& 0& 0& 1& 0\\0& 0&1& 0& 0& 0&1& 0& 1 \end{pmatrix}$,
  $Tr_1=(3,0,1), Tr_2=(3,0,2).$\ \
		\item
$As_2^4(3)= \begin{pmatrix}  1& 0& 0& 0& 1& -1& 1& 0& 0\\0&1& 0&1& 0& 1& 0& 1& 0\\0& 0&1& 0& 0& 1& 0& 0& 1  \end{pmatrix},$
$Tr_1=(2,0,2),\ Tr_2=(3,1,3)$.
		\item
 $As_2^5(3)= \begin{pmatrix} 1& 0& 1& 0& 1& 0& 0& -1& 0\\0& 1& 0& 1& 0& 1& 0& 1& 0\\0& 0& 0& 0& 0& 0& 1&1& 1  \end{pmatrix},$ $ Tr_1=(3,1,3), Tr_2=(2,0,2)$.
		\item
$As_2^6(3)=\begin{pmatrix}   1& 0& 0& 0& 0& 0& 0& 0& 0\\0& 0& 0& 0& 0& 0& 0& 1& 0\\0& 0& 0& 0& 0& 0& 0& 0& 1   \end{pmatrix},$ $Tr_1=(1,0,1), Tr_2=(1,0,2).$
		\item
 $As_2^7(3)=\begin{pmatrix}1& 0& 0& 0& 0& 0& 0& 0& 0\\0& 0& 0& 0& 0& 1& 0& 0& 0\\0& 0& 0& 0& 0& 0& 0& 0& 1\end{pmatrix},$ $Tr_1=(1,0,2), Tr_2=(1,0,1).$
		\item
$As_2^8(3)=\begin{pmatrix}  1& 0& 0& 0& 0& 0& 0& 0& 0\\0& 1& 0& 0& 0& 1& 0& 0& 0\\0& 0& 0& 0& 0& 0& 0& 0& 1  \end{pmatrix},$ $Tr_1=(1,0,2), Tr_2=(2,0,1).$
\end{enumerate}	
		\item Associative algebras with  $\lambda Tr_1=Tr_2,$ where $\lambda=3:$
\begin{enumerate}
		\item[(7)]
$As_{1,3}(3)=
		\begin{pmatrix}1& 0& 0& 0& 0& 0& 0& 0& 0\\0& 1& 0& 0& 0& 0& 0& 0& 0\\0& 0& 1& 0& 0& 0& 0& 0& 0\end{pmatrix},$
$Tr_1=(1,0,0), Tr_2=(3,0,0) $.
\end{enumerate}
	\item Associative algebras with  $\lambda Tr_1=Tr_2,$ where $\lambda=\frac{1}{3}:$
\begin{enumerate}
		\item[(8)]
$As_{1,1/3}(3)=
		\begin{pmatrix}1& 0& 0& 0& 0& 0& 0& 0& 0\\0& 0& 0&1& 0& 0& 0& 0& 0\\0& 0& 0& 0& 0& 0& 1& 0& 0\end{pmatrix},$
$Tr_1=(\frac{3}{2},0,0), Tr_2=(\frac{1}{2},0,0)$.
\end{enumerate}
	\item Associative algebras with  $\lambda Tr_1=Tr_2,$ where $\lambda=2:$
	\begin{enumerate}
		\item[(9)]
 $As_{1,2}(3)=\begin{pmatrix}  1& 0& 0& 0& 0& 0& 0& 0& 0\\0& 0& 0& 0& 0& 0& 0& 0& 0\\0& 0& 1& 0& 0& 0& 0& 0& 0   \end{pmatrix},\ \ Tr_1=(1,0,0), Tr_2=(2,0,0).$
	\end{enumerate}
	\item {Associative algebras with  $\lambda Tr_1=Tr_2,$ where $\lambda=\frac{1}{2}$:}
	\begin{enumerate}
		\item[(10)]
$As_{1,1/2}(3)=
		\begin{pmatrix}1& 0& 0& 0& 0& 0& 0& 0& 0\\0& 0& 0& 0& 0& 0& 0& 0& 0\\0& 0& 0& 0& 0& 0& 1& 0& 0\end{pmatrix},\ \ Tr_1=(2,0,0), Tr_2=(1,0,0);$
	\end{enumerate}
	\item {Associative algebras with  $\lambda Tr_1=Tr_2,$ where $\lambda=\frac{3}{2}:$}
	\begin{enumerate}
		\item[(11)]
   $As_{1, 3/2}(3)=
		\begin{pmatrix}1& 0& 0& 0&0& 0& 0& 0& 0\\0& 1& 0&1& 0& 0& 0& 0& 0\\0& 0&1& 0& 0&0& 0& 0& 0\end{pmatrix},$ $Tr_1=(2, 0, 0), \ \ Tr_2=(3, 0, 0).$
	\end{enumerate}
	\item {Associative algebras with $\lambda Tr_1=Tr_2,$ where $\lambda=\frac{2}{3}:$}
	\begin{enumerate}
		\item[(12)]
$As_{1,2/3}(3)=
		\begin{pmatrix}1& 0& 0& 0& 0& 0& 0& 0& 0\\0&1& 0&1& 0& 0& 0& 0& 0\\0& 0& 0& 0& 0& 0&1& 0& 0\end{pmatrix},$
 $Tr_1=(3,0,0), Tr_2=(2,0,0)$.
	\end{enumerate}
\item	Associative algebras with $\lambda Tr_1=Tr_2,$ where $\lambda=1:$
	\begin{enumerate}
		\item[(13)]
$As_{1,1}^{1}(3)(t)=\begin{pmatrix}1& 0& 0& 0& t& 0& 0& 0& 0\\0&1& 0& 1& 0& 0& 0& 0& 0\\0& 0& 0& 0& 0& 0& 0& 0& 0   \end{pmatrix}$ $(As_{1,1}^{1}(3)(t)\cong As_{1,1}^{1}(3)(a^2t),\ a\in \mathbb{F}^*), \  t \in \mathbb{F},$ \ $Tr_1=Tr_2=(2,0,0).$
		\item[(14)]
 $As_{1,1}^{2}(3)(t)=\begin{pmatrix}1& 0& 0& 0& t& 0& 0& 0& 0\\0&1& 0&1& 0& 0& 0& 0& 0\\0& 0& 0& 0& 0& 0& 0& 0& 1   \end{pmatrix}$ $(As_{1,1}^{2}(3)(t)\cong As_{1,1}^{2}(3)(a^2t),\ a\in \mathbb{F}^*), \ t \in \mathbb{F},$ \ $Tr_1=Tr_2=(2,0,1).$
		\item[(15)]
  $As_{1,1}^{3}(3)=\begin{pmatrix}1& 0& 0& 0& 0& 0& 0& 0& 0\\0&1& 0& 1& 0& 0& 0& 0& 1\\0& 0&1& 0& 0& 0&1& 0& 0  \end{pmatrix},$\ $Tr_1=Tr_2=(3,0,0).$ 
		\item[(16)]
  $As_{1,1}^{5}(3)=\begin{pmatrix}1& 0& 0& 0& 0& 0& 0& 0& 0\\0&1& 0& 1& 0& 0& 0& 0& 0\\0& 0&1& 0& 0& 0&1& 0& 0  \end{pmatrix},\  Tr_1=Tr_2=(3,0,0).$
  	
 \item[(17)] $As_{1,1}^{6}(3)(-1)=\begin{pmatrix} 1& 0& 0& 0& 0& 0& 0& 0&-1\\0& 1& 0& 1& 0& 1& 0& 1& 1\\0& 0& 1& 0& 0& 0&1& 0&2 \end{pmatrix},$
 $Tr_1= Tr_2=(3,0,3).$
 \item[(18)]
$As_{1,1}^{7}(3)(-1)=\begin{pmatrix} 1& 0& 0& 0& 0& 0& 0& 0&-1\\0& 1& 0& 1& 0& 1& 0& 1& 0\\0& 0& 1& 0& 0& 0&1& 0&2 \end{pmatrix},$
$Tr_1= Tr_2=(3,0,3).$\\
		\item[(19)]
$As_{1,1}^{8}(3)(t)=\begin{pmatrix}1& 0& 0& 0& 1& 0& 0& 0& t\\0& 1& 0& 1& 0& 0& 0& 0& t\\0& 0&1& 0& 0& 1&1& 1& 1\end{pmatrix},$ $ t \in \mathbb{F}.$

Note that $As_{1,1}^{8}(3)(t)\simeq As_{1,1}^8(3)(s)$, whenever $8s+1\neq 0$, there exists $a\in \mathbb{F}$
such that $a^2=1+8(t(8s+1)+s)$ and $1-4(8t+1+a)\neq 0$.
		\item[(20)]  
 $As_{1,1}^{9}(3)(t)= \begin{pmatrix}1& 0& 0& 0& 1& 0& 0& 0& t\\0&1& 0&1& 0& 0& 0& 0& t\\0& 0&1& 0& 0&1&1& 1& 0\end{pmatrix}$ $(As_{1,1}^{9}(3)(t)\cong As_{1,1}^{9}(3)(a^2t), \ a\in \mathbb{F}^*, \ t \in \mathbb{F},$  \ $Tr_1= Tr_2=(3,1,0).$
		\item[(21)] $As_{1,1}^{10}(3)=\begin{pmatrix}  1& 0& 0& 0& 0& 0& 0& 0& 0\\0& 0& 0& 0& 0& 0& 0& 0& 0\\0& 0& 0& 0& 0& 0& 0& 0& 0   \end{pmatrix},$
		\ $Tr_1= Tr_2=(1,0,0).$
		\item[(22)] $As_{1,1}^{11}(3)=\begin{pmatrix}1& 0& 0& 0& 0& 0& 0& 0& 0\\0& 0& 0& 0& 0& 1& 0& 1& 0\\0& 0& 0& 0& 0& 0& 0& 0& 1\end{pmatrix},\ \ Tr_1= Tr_2=(1,0,2).$
		\item[(23)] $As_{1,1}^{12}(3)=
		\begin{pmatrix}1& 0& 0& 0& 0& 0& 0& 0& 0\\0& 1& 0& 0& 0& 0& 0& 0& 0\\0& 0& 0& 0& 0& 0& 1& 0& 0\end{pmatrix},$\ \
  $Tr_1= Tr_2=(2,0,0) $.
		\item[(24)] $As_{1,1}^{13}(3)=\begin{pmatrix}  0& 0& 0& 0& 0& 0& 0& 0& 0\\1& 0& 0& 0& 0& 0& 0& 0& 0\\0& 0& 0& 0& 0& 0& 0& 0& 1   \end{pmatrix},$ \ \ $ Tr_1= Tr_2=(0,0, 1).$ \\
		\item[(25)]
		$As_{1,1}^{14}(3)=\begin{pmatrix}0&0&0&0&0&0&
		0&0&1\\ 0&0&0&0&0&1&0&1&t\\ 0&0&0&0&0&0&0&0&1\end{pmatrix}$, $t\in \{0,1\}, \ Tr_1=  Tr_2=(0,0, 2).$
\item[(26)] $As_{1,1}^{17}(3)(t)=\begin{pmatrix}0&0&t&0&t&0&t&0&0\\
    1&0&0&0&0&t&0&t&0\\
    0&1&0&1&0&0&0&0&t\end{pmatrix}$ $(As_{1,1}^{17}(3)(t)\cong As_{1,1}^{17}(3)(a^3t),\ a\in \mathbb{F}^*)$, provided that  $u^3=t$ has no solution in $\mathbb{F},\ t\in \mathbb{F}$,
    $ Tr_1=  Tr_2=(0,0, 3t).$	
 	\end{enumerate}
		\item {Associative algebras with $Tr_1=Tr_2=0$:}
	\begin{enumerate}
		\item[(27)]
$As_0^{1}(3)(1)= \begin{pmatrix}0& 0& 0& 0& 0&1& 0& 1& 0\\0& 0& 0& 0& 0& 0& 0& 0& 0\\0& 0& 0& 0& 0& 0& 0& 0& 0\end{pmatrix}$\\
 if $u^2=-1$ has a solution in $\mathbb{F}$ then $As_0^{1}(3)(1)\cong As_0^{4}(3)(1)$.
\item[(28)] $As_0^{1}(3)(-1)= \begin{pmatrix}0& 0& 0& 0& 0&-1& 0& 1& 0\\0& 0& 0& 0& 0& 0& 0& 0& 0\\0& 0& 0& 0& 0& 0& 0& 0& 0\end{pmatrix}.$
\item[(29)]
$As_0^{2}(3)=\begin{pmatrix}    0& 0& 0& 0& 0& 0& 0& 0& 0\\1& 0& 1& 0& 0& 0& -1& 0&0 \\0& 0& 0& 0& 0& 0& 0& 0& 0  \end{pmatrix}.$
 	\item[(30)]
 	$As_0^{3}(3)=\begin{pmatrix}0&0&0&0&0&1&
 		0&1&0\\ 0&0&0&0&0&0&0&0&1\\ 0&0&0&0&0&0&0&0&0\end{pmatrix}$.
		\item[(31)]
$As_0^{4}(3)(t)=\begin{pmatrix}  0& 0& 0& 0& 0& 0& 0& 0& 0\\1& 0& 0& 0& 0& 0& 0& 0& t\\0& 0& 0& 0& 0& 0& 0& 0& 0   \end{pmatrix}$ $(As_0^{4}(3)(t)\cong As_0^{4}(3)(a^2t),\  a\in \mathbb{F}^*),$ $ t\in \mathbb{F}.$
		\item[(32)] 
$As_0^{5}(3)(t)=\begin{pmatrix}    0& 0& 0& 0& 0& 0& 0& 0& 0\\1& 0& 0& 0& 0& 0& 1& 0&t \\0& 0& 0& 0& 0& 0& 0& 0& 0  \end{pmatrix},$\\
if $u^2=-1$ has a solution in $\mathbb{F}$, then $As_0^{5}(3)(t)\cong As_0^{5}(3)(-t)), \ t\in \mathbb{F}$.
	\end{enumerate}	
\end{itemize}
\end{theor}

\begin{proof} Here is a procedure to run for the proof.

Let $\mathbb{A}$ be a two-dimensional associative algebra and \[As(2)=\left(\begin{array}{cccc} \alpha_1 & \alpha_2 & \alpha_3 &\alpha_4\\ \beta_1 & \beta_2 & \beta_3 &\beta_4\end{array}\right)\]be its MSC with respect to the basis $\mathbf{e}=\left(\mathrm{e}_1, \mathrm{e}_2\right)$.

Then a three-dimensional algebra containing $As(2)$ as a subalgebra, with respect to the extended basis $\mathbf{e}=\left(\mathrm{e}_1, \mathrm{e}_2, \mathrm{e}_3\right)$ can be written as follows:
$$As(3)=\begin{pmatrix}\alpha_1& \alpha_2& x_3& \alpha_3& \alpha_4& x_6& x_7& x_8& x_9\\\beta_1& \beta_2& y_3&\beta_3& \beta_4& y_6& y_7& y_8& y_9\\0& 0& z_3& 0& 0& z_6&z_7& z_8& z_9\end{pmatrix}. $$
Now we impose the condition (\ref{AA2}) on $As(3)$, which yields a system of equations in the variables $x_i, y_j, z_k$ (where $i, j, k=3,6,... ,9$) with the coefficients $\alpha_s, \beta_t,$ ($s, t=1, 2, 3, 4$) regarded as known. As the two-dimensional subalgebra $As(2)$ we consider both the trivial algebra and the algebras from Theorem \ref{Ass}.
So, our procedure to classify three-dimensional associative algebras consists of the following steps:
\begin{itemize}
\item Solve the system of equations for $x_i, y_j, z_k$, where $i, j, k=3,6,... ,9$ to obtain a set of three-dimensional algebras containing the chosen two-dimensional subalgebra;
\item Check for isomorphisms within each group and eliminate redundancies, if any;
\item Classify three dimensional associative algebras which have no two-dimensional subalgebras, but have one-dimensional subalgebras;
\item Classify three-dimensional associative algebras which have only trivial subalgebras;
\item Compile the final list of three-dimensional associative algebras.
\end{itemize}
\end{proof}
\section{List of three-dimensional complex associative algebras}
In this section, we present the list of representatives of the isomorphism classes of three-dimensional complex associative algebras from \cite{KSTT}, along with their trace vectors. Note that the list of \cite{KSTT} agrees with that of  \cite{FP}.
The authors divide three-dimensional associative algebras into several groups: unital and non-unital, and the non-unital algebras, in their turn, are divided into three types: curled, waved and straight algebras. According to them an algebra
$\mathbb{A}$ is said to be \emph{curled} if $\mathrm{x}$ and $\mathrm{x}^2$ are linearly dependent for any $\mathrm{x} \in \mathbb{A}$, it is called \emph{waved} if it is not curled but $\mathrm{x},\ \mathrm{x}^2$ and $\mathrm{x}^3$ are linearly dependent for any $\mathrm{x} \in \mathbb{A}$, and it is
\emph{straight} otherwise. We keep the authors notations for the algebras as they are.

	 \begin{itemize}
 \item Unital algebras:
 \begin{enumerate}
  \item $U_0^3=\begin{pmatrix}1&0&0&0&0&0&0&0&0\\ 0&1&0&1&0&0&0&0&0\\ 0&0&1&0&0&0&1&0&0\end{pmatrix}, \ Tr_1=(3, 0, 0), \ \ Tr_2=(3, 0, 0).$
\item $U_1^3=\begin{pmatrix}1&0&0&0&0&0&0&0&1\\ 0&1&0&1&0&1&0&-1&0\\ 0&0&1&0&0&0&1&0&0\end{pmatrix}, \ Tr_1=(3, 0, 1), \ \ Tr_2=(3, 0, -1).$
\item $U_2^3=\begin{pmatrix}1&0&0&0&0&0&0&0&0\\ 0&0&0&0&1&0&0&0&0\\ 0&0&0&0&0&0&0&0&1\end{pmatrix}, \ Tr_1=(1, 1, 1), \ \ Tr_2=(1, 1, 1).$
\item $U_3^3=\begin{pmatrix}1&0&0&0&0&0&0&0&0\\ 0&0&0&0&1&0&0&0&0\\ 0&0&0&0&0&1&0&1&0\end{pmatrix}, \ Tr_1=(1, 2, 0), \ \ Tr_2=(1, 2, 0).$
\item $U_4^3=\begin{pmatrix}1&0&0&0&0&0&0&0&0\\ 0&1&0&1&0&0&0&0&0\\ 0&0&1&0&1&0&1&0&0\end{pmatrix}, \ Tr_1=(3, 0, 0), \ \ Tr_2=(3, 0, 0).$
  \end{enumerate}
  \item Curled algebras:
  \begin{enumerate}
  \item[(6)] The trivial algebra is denoted by $C_0^3$.
  \item[(7)] $C_1^3=\begin{pmatrix}0&0&0&0&0&1&0&-1&0\\ 0&0&0&0&0&0&0&0&0\\ 0&0&0&0&0&0&0&0&0\end{pmatrix}, \ Tr_1=(0, 0, 0), \ \ Tr_2=(0, 0, 0).$
  \item[(8)] $C_2^3=\begin{pmatrix}0&0&0&1&0&0&0&0&0\\ 0&0&0&0&1&0&0&0&0\\ 0&0&0&0&0&0&0&1&0\end{pmatrix}, \ Tr_1=(0, 2, 0), \ \ Tr_2=(0, 2, 0).$
  \item[(9)] $C_3^3=\begin{pmatrix}0&0&0&0&0&0&1&0&0\\ 0&0&0&0&0&0&0&1&0\\ 0&0&0&0&0&0&0&0&1\end{pmatrix}, \ Tr_1=(0, 0, 1), \ \ Tr_2=(0, 0, 3).$
  \item[(10)] $C_4^3=\begin{pmatrix}0&0&1&0&0&0&0&0&0\\ 0&0&0&0&0&1&0&0&0\\ 0&0&0&0&0&0&0&0&1\end{pmatrix}, \ Tr_1=(0, 0, 3), \ \ Tr_2=(0, 0, 1).$
  \end{enumerate}
  \item Straight algebras:
   \begin{enumerate}
   \item[(11)] $S_1^3=\begin{pmatrix}0&0&0&0&0&0&0&0&0\\ 1&0&0&0&0&0&0&0&0\\ 0&1&0&1&0&0&0&0&0\end{pmatrix}, \ Tr_1=(0, 0, 0), \ \ Tr_2=(0, 0, 0).$
  \item[(12)] $S_2^3=\begin{pmatrix}1&0&0&0&0&0&0&0&0\\ 0&0&0&0&0&0&0&0&0\\ 0&0&0&0&1&0&0&0&0\end{pmatrix}, \ Tr_1=(1, 0, 0), \ \ Tr_2=(1, 0, 0).$
  \item[(13)] $S_3^3=\begin{pmatrix}1&0&0&0&0&0&0&0&0\\ 0&0&0&0&1&0&0&0&0\\ 0&0&0&0&0&0&0&0&0\end{pmatrix}, \ Tr_1=(1, 1, 0), \ \ Tr_2=(1, 1, 0).$
  \item[(14)] $S_4^3=\begin{pmatrix}1&0&0&0&0&0&0&0&0\\ 0&1&0&1&0&0&0&0&0\\ 0&0&0&0&0&0&0&0&0\end{pmatrix}, \ Tr_1=(2, 0, 0), \ \ Tr_2=(2, 0, 0).$
  \end{enumerate}
  \item Waved algebras:
  \begin{enumerate}
   \item[(15)] $W_1^3=\begin{pmatrix}0&0&0&0&0&0&0&0&1\\ 0&0&0&0&0&0&0&0&0\\ 0&0&0&0&0&0&0&0&0\end{pmatrix}, \ Tr_1=(0, 0, 0), \ \ Tr_2=(0, 0, 0).$
  \item[(16)] $W_2^3=\begin{pmatrix}0&0&0&0&0&0&0&1&0\\ 0&0&0&0&0&0&0&0&0\\ 0&0&0&0&0&0&0&0&0\end{pmatrix}, \ Tr_1=(0, 0, 0), \ \ Tr_2=(0, 0, 0).$
  \item[(17)] $W_3^3=\begin{pmatrix}0&0&0&0&1&0&0&k&1\\ 0&0&0&0&0&0&0&0&0\\ 0&0&0&0&0&0&0&0&0\end{pmatrix}, \ k\in\{x+iy | x>0 \ \mbox{or}\ x=0, y\geq0\},$

      \hfill $Tr_1=(0, 0, 0), \ \ Tr_2=(0, 0, 0).$
  \item[(18)] $W_4^3=\begin{pmatrix}1&0&0&0&0&0&0&0&0\\ 0&0&0&0&0&0&0&0&0\\ 0&0&0&0&0&0&0&0&0\end{pmatrix}, \ Tr_1=(1, 0, 0), \ \ Tr_2=(1, 0, 0).$
  \item[(19)] $W_5^3=\begin{pmatrix}0&0&0&0&0&0&0&0&0\\ 0&0&0&0&0&0&0&1&0\\ 0&0&0&0&0&0&0&0&1\end{pmatrix}, \ Tr_1=(0, 0, 1), \ \ Tr_2=(0, 0, 2).$
  \item[(20)] $W_6^3=\begin{pmatrix}0&0&0&0&0&0&0&0&0\\ 0&0&0&0&0&1&0&0&0\\ 0&0&0&0&0&0&0&0&1\end{pmatrix}, \ Tr_1=(0, 0, 2), \ \ Tr_2=(0, 0, 1).$
  \item[(21)] $W_7^3=\begin{pmatrix}1&0&0&0&0&0&0&0&0\\ 0&0&0&0&0&0&0&1&0\\ 0&0&0&0&0&0&0&0&1\end{pmatrix}, \ Tr_1=(1, 0, 1), \ \ Tr_2=(1, 0, 2).$
  \item[(22)] $W_8^3=\begin{pmatrix}1&0&0&0&0&0&0&0&0\\ 0&0&0&0&0&1&0&0&0\\ 0&0&0&0&0&0&0&0&1\end{pmatrix}, \ Tr_1=(1, 0, 2), \ \ Tr_2=(1, 0, 1).$
  \item[(23)] $W_9^3=\begin{pmatrix}0&1&0&1&0&0&0&0&0\\ 0&0&0&0&1&0&0&0&0\\ 0&0&0&0&0&0&0&1&0\end{pmatrix}, \ Tr_1=(0, 3, 0), \ \ Tr_2=(0, 2, 0).$
  \item[(24)] $W_{10}^3=\begin{pmatrix}0&1&0&1&0&0&0&0&0\\ 0&0&0&0&1&0&0&0&0\\ 0&0&0&0&0&1&0&0&0\end{pmatrix}, \ Tr_1=(0, 2, 0), \ \ Tr_2=(0, 3, 0).$
  \end{enumerate}
  \end{itemize}
\begin{itemize}
\item Comparison with the list of complex associative algebras in \cite{KSTT}.
\begin{center}	
{\footnotesize{\begin{tabular}{|c|c|c|c|c|c|c|c|c}
 \hline
Algebra            & {Algebra in}               &Base change $g$& Isomorphism\\
from \cite{KSTT}&   this paper           &&  \\ \hline
$U_0^3$   & $As_{1,1}^{5}(3)$                 &=& $U_0^3=As_{1,1}^{5}(3)$\\
 $ U_1^3$ &  $ As_2^2(3)\left(0\right)$&$ \begin{pmatrix}1&0&1\\ 0&-2&-2\\ 0&0&-2\end{pmatrix}$& $U_1^3=g^{-1}As_2^2(3)\left(0\right)g^{\otimes 2}$  \\
 $ U_2^3$  &     $As_{1,1}^{2}(3)(1)$     &$ \begin{pmatrix}\frac{1}{2}&\frac{1}{2}&0\\ \frac{1}{2}&-\frac{1}{2}&0\\ 0&0&1\end{pmatrix}$& $U_2^3=g^{-1}As_{1,1}^{2}(3)(1)g^{\otimes 2}$  \\
 $ U_3^3$ & $As_{1,1}^{11}(3)$           &$ \begin{pmatrix}1&0&1\\ 1&0&-1\\ 0&1&0\end{pmatrix}$& $U_3^3=g^{-1}As_{1,1}^{11}(3)g^{\otimes 2}$ \\
 $ U_4^3$ &    $As_{1,1}^{3}(3)$       &$ \begin{pmatrix}1&0&0\\ 0&1&1\\ 0&1&0\end{pmatrix}$&$U_4^3=g^{-1}As_{1,1}^{3}(3)g^{\otimes 2}$ \\
$ C_1^3$  & $As_0^{1}(3)(-1)$               &$ \begin{pmatrix}-1&0&0\\ 0&1&0\\ 0&0&1\end{pmatrix}$& $ C_1^3=g^{-1}As_0^{1}(3)(-1)g^{\otimes 2}$\\
$ C_2^3$  &$As_{1,1}^{12}(3)$         &$ \begin{pmatrix}0&1&0\\ 1&0&0\\ 0&0&1\end{pmatrix}$& $C_2^3=g^{-1}As_{1,1}^{12}(3)g^{\otimes 2}$ \\
$ C_3^3$  &$ As_{1,3}(3)$         &$\begin{pmatrix}0&0&1\\ 1&0&0\\ 0&1&0\end{pmatrix}$& $C_3^3=g^{-1}As_{1,3}(3)g^{\otimes 2}$  \\
$ C_4^3$  &$ As_{1,1/3}(3)$          &$ \begin{pmatrix}0&0&1\\ 0&1&0\\ 1&0&0\end{pmatrix}$& $C_4^3=g^{-1}As_{1,1/3}(3)g^{\otimes 2}$   \\
$ S_1^3$  & $As_0^{3}(3)(1)$               &$ \begin{pmatrix}0&0&1\\ 0&1&0\\ 1&0&0\end{pmatrix}$& $S_1^3=g^{-1}As_0^{3}(3)(1)g^{\otimes 2}$   \\
 $ S_2^3$  &$As_{1,1}^{13}(3)$ &$ \begin{pmatrix}0&1&0\\ 0&0&1\\ 1&0&0\end{pmatrix}$& $S_2^3=g^{-1}As_{1,1}^{13}(3)g^{\otimes 2}$   \\
$ S_3^3$  &$As_{1,1}^{1}(3)(1)$        &$\begin{pmatrix}\frac{1}{2}&\frac{1}{2}&0\\ \frac{1}{2}&-\frac{1}{2}&0\\ 0&0&1\end{pmatrix}$& $S_3^3=g^{-1}As_{1,1}^{1}(3)(1)g^{\otimes 2}$   \\
$S_4^3$  &$As_{1,1}^{1}(3)(0)$         &=& $S_4^3=As_{1,1}^{1}(3)(0)$ \\
 $W_1^3$ & $As_0^{4}(3)(0)$           &$ \begin{pmatrix}0&0&1\\ 1&0&0\\ 0&1&0\end{pmatrix}$& $W_1^3=g^{-1}As_0^{4}(3)(0)g^{\otimes 2}$   \\
$W_2^3$  & $As_0^5(3)(0)$               &$ \begin{pmatrix}0&-1&0\\ -1&0&0\\ 0&1&1\end{pmatrix}$& $W_2^3\cong As_0^5(3)(0)$  \\
   \hline
   \end{tabular}}}
\end{center}

\begin{center}	
{\footnotesize{\begin{tabular}{|c|c|c|c|c|c|c|c|c}
 \hline
Algebra            & {Algebra in}               &Base change $g$& Isomorphism\\
from \cite{KSTT}&  this paper           &&  \\ \hline
 $W_3^3(k), k\neq0$  &$As_0^{5}(3)\left(\frac{1}{k^2}\right)$           &$ \begin{pmatrix}0&0&-\frac{1}{k}\\ \frac{1}{k^2}&0&1\\ 0&1&0\end{pmatrix}$& $W_3^3(k)=g^{-1}As_0^{5}(3)\left(\frac{1}{k^2}\right)g^{\otimes 2}$   \\
 $W_3^3(0)$  &$As_0^{4}(3)(1)$           &$ \begin{pmatrix}0&0&1\\ 1&0&0\\ 0&1&0\end{pmatrix}$& $W_3^3(0)=g^{-1}As_0^{4}(3)(1)g^{\otimes 2}$   \\
$W_4^3$ &$As_{1,1}^{10}(3)$           & = & $W_4^3=As_{1,1}^{10}(3)$   \\
$W_5^3$  &$As_{1,2}(3)$        &$\begin{pmatrix}0&0&1\\ 1&0&0\\ 0&1&0\end{pmatrix}$& $W_5^3=g^{-1} As_{1,2}(3)g^{\otimes 2}$   \\
 $W_6^3$  &$As_{1,1/2}(3)$        &$\begin{pmatrix}0&0&1\\ 1&0&0\\ 0&1&1\end{pmatrix}$& $W_6^3=g^{-1}As_{1,1/2}(3)g^{\otimes 2}$   \\
$W_7^3$ &$ As_2^6(3)$        &= & $W_7^3= As_2^6(3)$   \\
$W_8^3$  &$ As_2^7(3)$        & =& $W_8^3=As_2^7(3)$  \\
$W_9^3$ &$ As_{1,2/3}(3)$        &$\begin{pmatrix}0&1&0\\ 1&0&0\\ 0&0&1\end{pmatrix}$& $W_9^3=g^{-1}As_{1,2/3}(3)g^{\otimes 2}$   \\
$W_{10}^3$  &$As_{1,3/2}(3)$        &$\begin{pmatrix}0&1&0\\ 1&0&0\\ 0&0&1\end{pmatrix}$& $W_{10}^3=g^{-1}As_{1,3/2}(3)g^{\otimes 2}$   \\
   \hline
\end{tabular}}}
\end{center}
\emph{}
\\
\item \textbf{Conclusions}.
Here is a list of three-dimensional complex associative algebras to be added to the list of \cite{KSTT} to have a complete list of such algebras:
\begin{itemize}
\item \textbf{Unital}:
	 $\bullet\ As_{1,1}^{6}(3)(-1);$\ \ $\bullet\ As_{1,1}^{7}(3)(-1);$
$\bullet\ As_{1,1}^{9}(3)(1)$.

In the case $\mathbb{F}=\mathbb{C}$, due to $As_{1,1}^8(3)(s)\simeq As_{1,1}^8(3)(t),$ whenever $8t+1\neq 0$, there exists $a$ such that $a^2=1+8(s(8t+1)+t)$ and $1-4(8t+1+a)\neq 0$, setting $s=0$ we have $As_{1,1}^8(3)(0)\simeq As_{1,1}^8(3)(t)$, where $8t+1\neq 0,$  $a^2= 8t+1$ and
$1-4(a^2+ a)\neq 0$. Therefore, in this case we get the following two non-isomorphic algebras: \\ $\bullet\ As_{1,1}^8(3)(0)$;\ \ $\bullet\ As_{1,1}^8(3)(-\frac{1}{8})$.
\item \textbf{Waved}: 
$\bullet\ As_2^4(3);$ \ \ $\bullet\ As_2^5(3)$;
\ \ $\bullet\ As_0^{2}(3)$.
\item \textbf{Straight}:
  $\bullet\ As_2^{8}(3);$\ \ $\bullet\ As_{1,1}^{14}(3)(t),\ t \in \{0,1\}$.

\vspace{0.5cm}

Here is the list of arguments that can be used to distinguish the extra algebras:
\begin{itemize}
  		\item $M_2(n)=\left\{A \in M_{n\times n^2}(\mathbb{F}):\ \  \{Tr_1(A), \ Tr_2(A)\} \ \mbox{is linearly independent} \right\};$
  		\item $M_{1,\lambda}(n)=\left\{A \in M_{n\times n^2}(\mathbb{F}):\ \  \lambda Tr_1(A)= Tr_2(A) \ \mbox{and non of}\  Tr_1(A), Tr_2(A) \ \mbox{is zero}\right\};$
  		\item $M_{1,0}(n)=\left\{A \in M_{n\times n^2}(\mathbb{F}): \  \ Tr_1(A)\ \mbox{is nonzero and}\ Tr_2(A)\ \mbox{is zero }\right\};$
  		\item $M_{0}(n)=\left\{ A \in M_{n\times n^2}(\mathbb{F}): \ \mbox{the both}\ Tr_1(A), Tr_2(A) \ \mbox{are zero} \right\}.$
  	\end{itemize}

\begin{enumerate}
\item $i$ - the number of the left ideals  (designated as $(1): i$);
  \item  $(a,b),$ where $a$ and $b$ are the numbers of one-dimensional and two-dimensional left ideals, respectively (designated as $(2): [a,b]$);
  \item  $j$ - the number of the right ideals  (designated as $(3): j$);
\item Decomposability (designated as (4): D/ID);
\item The number of idempotents (designated as (5): NI);
\end{enumerate}
{\tiny{
\begin{table}[h]
\caption{The comparison of the extra unital algebras.}
\begin{center}	
\begin{tabular}{|c|c|c|c|c|c|}
\hline
\diagbox{\cite{KSTT} }{\emph{} }  & $As_{1,1}^{6}(3)(-1)$  & $ As_{1,1}^{7}(3)(-1)$                &$ As_{1,1}^{9}(3)(1)$ &$ As_{1,1}^{8}(3)(0)$ & $ As_{1,1}^{8}(3)\left(-\frac{1}{8}\right)$\\ \hline
 $U_0^3$             &\diagbox{(4): D}{\emph{} (4): ID \emph{}}  &     \diagbox{(4): D}{\emph{} (4): ID \emph{}}                 &    \diagbox{(4): D}{\emph{} (4): ID \emph{}}      &     \diagbox{(4): D}{\emph{} (4): ID \emph{}}&\diagbox{(4): D}{\emph{} (4): ID } \\ \hline
$U_1^3$ &\diagbox{$M_{1,1}(3)$}{\emph{} $M_0(3)$ \emph{}} &  \diagbox{$(1):4$}{$(1):8$}     &   \diagbox{$(3):7$}{$(3):9$}&    \diagbox{$(2):[1,1]$}{$(2):[2,1]$} &\diagbox{(5): NI=many}{\emph{} (5): NI=4 } \\ \hline
$U_2^3$ &\diagbox{$M_{1,1}(3)$}{\emph{} $M_0(3)$ \emph{}} &  \diagbox{$(1):4$}{$(1):8$}     &   \diagbox{$(3):7$}{$(3):9$}&    \diagbox{$(2):[1,1]$}{$(2):[2,1]$} & \diagbox{(4): D}{\emph{} (4): ID \emph{}}\\ \hline
$U_3^3$ & \diagbox{$M_{1,1}(3)$}{\emph{} $M_0(3)$ \emph{}} &    \diagbox{$(1):4$}{$(1):8$}    &   \diagbox{$(3):7$}{$(3):9$} &    \diagbox{$(2):[1,1]$}{$(2):[2,1]$} &\diagbox{(4): D}{\emph{} (4): ID \emph{}} \\ \hline
$U_4^3$ & \diagbox{$M_{1,1}(3)$}{\emph{} $M_0(3)$ \emph{}} &    \diagbox{$(1):4$}{$(1):8$}    &   \diagbox{$(3):7$}{$(3):9$} &    \diagbox{$(2):[1,1]$}{$(2):[2,1]$} &\diagbox{(4): D}{\emph{} (4): ID \emph{}} \\ \hline
\end{tabular}
\end{center}
\end{table}}}

\begin{table}[h]
\caption{The comparison of the extra waved algebras.}
{\tiny{\begin{tabular}{|c|c|c|c|c|}
\hline
\diagbox{\cite{KSTT} }{\emph{} }  & $As_2^4(3)$                                    & $As_2^5(3)$   & $As_0^2(3)$  \\ \hline
        $W_1^3$                   &\diagbox{$M_2(3)$}{\emph{} $M_0(3)$ \emph{}}    &     \diagbox{$M_2(3)$ }{$M_{1,1}(3)$} & \diagbox{$M_2(3)$}{\emph{} $M_0(3)$ \emph{}}   \\ \hline
 $W_2^3$                          &\diagbox{$M_{1,1}(3)$}{\emph{} $M_0(3)$ \emph{}} &  \diagbox{$(1):4$}{$(1):8$} & \diagbox{$M_{1,1}(3)$}{\emph{} $M_0(3)$ \emph{}}  \\ \hline
$W_3^3$                           &\diagbox{$M_{1,1}(3)$}{\emph{} $M_0(3)$ \emph{}} &  \diagbox{$(1):4$}{$(1):8$}  &  \diagbox{$M_{1,1}(3)$}{\emph{} $M_0(3)$ \emph{}}  \\ \hline
\qquad $W_4^3$    \qquad\qquad                       & \diagbox{$M_{1,1}(3)$}{\emph{} $M_0(3)$ \emph{}} &    \diagbox{$(1):4$}{$(1):8$}&\diagbox{$M_{1,1}(3)$}{\emph{} $M_0(3)$ \emph{}}  \\ \hline
        $W_5^3$                   &\diagbox{$M_2(3)$}{\emph{} $M_0(3)$ \emph{}}    &     \diagbox{$M_2(3)$ }{$M_{1,1}(3)$} &  \diagbox{$M_{2}(3)$}{\emph{} $M_0(3)$ \emph{}}   \\ \hline
 $W_6^3$                          &\diagbox{$M_{1,1}(3)$}{\emph{} $M_0(3)$ \emph{}} &  \diagbox{$(1):4$}{$(1):8$} & \diagbox{$M_{1,1}(3)$}{\emph{} $M_0(3)$ \emph{}}  \\ \hline
$W_7^3$                           &\diagbox{$M_{1,1}(3)$}{\emph{} $M_0(3)$ \emph{}} &  \diagbox{$(1):4$}{$(1):8$}  &  \diagbox{$M_{1,1}(3)$}{\emph{} $M_0(3)$ \emph{}}  \\ \hline
$W_8^3$                           & \diagbox{$M_{1,1}(3)$}{\emph{} $M_0(3)$ \emph{}} &    \diagbox{$(1):4$}{$(1):8$}& \diagbox{$M_{1,1}(3)$}{\emph{} $M_0(3)$ \emph{}}\\ \hline
        $W_9^3$                   &\diagbox{$M_2(3)$}{\emph{} $M_0(3)$ \emph{}}    &     \diagbox{$M_2(3)$ }{$M_{1,1}(3)$} & \diagbox{$M_{2}(3)$}{\emph{} $M_0(3)$ \emph{}}  \\ \hline
 $W_{10}^3$                          &\diagbox{$M_{1,1}(3)$}{\emph{} $M_0(3)$ \emph{}} &  \diagbox{$(1):4$}{$(1):8$} & \diagbox{$M_{1,1}(3)$}{\emph{} $M_0(3)$ \emph{}} \\ \hline
\end{tabular}}}
\end{table}

\begin{table}[h]
\caption{Comparison of the extra straight algebras.}
{\tiny{\begin{tabular}{|c|c|c|c|c|}
\hline
\diagbox{\cite{KSTT} }{\emph{} }  & $As_2^8(3)$   & $As_{1,1}^{14}(3)(0)$              & $As_{1,1}^{14}(3)(1)$ \\ \hline
$S_1^3$             &\diagbox{\emph{} $M_0(3)$ \emph{}}{$M_2(3)$}  &     \diagbox{$M_{0}(3)$ }{$M_{1,1}(3)$}                  &    \diagbox{$M_{0}(3)$ }{$M_{1,1}(3)$}        \\ \hline
$S_2^3$  &\diagbox{\emph{} $M_{1,1}(3)$ \emph{}}{$M_{2}(3)$} &  \diagbox{$(1):8$}{$(1):4$}     &   \diagbox{$(3):9$}{$(3):7$}\\ \hline
 $S_3^3$ &\diagbox{\emph{} $M_{1,1}(3)$ \emph{}}{$M_{2}(3)$} &  \diagbox{$(1):8$}{$(1):4$}     &   \diagbox{$(3):9$}{$(3):7$}\\ \hline
$S_4^3$ & \diagbox{\emph{} $M_{1,1}(3)$ \emph{}}{$M_{2}(3)$} &    \diagbox{$(1):8$}{$(1):4$}    &   \diagbox{$(3):9$}{$(3):7$} \\ \hline
\end{tabular}}}
\end{table}

\newpage

\section{List of three-dimensional nilpotent associative algebras}
The list of three-dimensional nilpotent associative algebras over any base field was given in  \cite{GR} as follows:
\begin{enumerate}
 \item $A_{3,2}=\begin{pmatrix}0&0&0&0&0&0&0&0&0\\ 1&0&0&1&0&0&0&0&0\\ 0&0&0&0&0&0&1&0&0\end{pmatrix}.$
\item $A_{3,3}^\delta=\begin{pmatrix}0&0&0&0&0&0&0&0&1\\ 0&0&0&0&0&0&0&0&0\\ 1&0&0&0&\delta&0&0&0&0\end{pmatrix},$
$(A_{3,3}^\delta \cong A_{3,3}^{a^2\delta}, \ a, \delta \in \mathbb{F}^*).$
\item $A_{3,4}^\delta=\begin{pmatrix}0&0&0&0&0&0&0&0&0\\ 0&0&0&0&0&0&0&0&0\\ 1&1&0&0&\delta&0&0&0&0\end{pmatrix}, \ \delta \in \mathbb{F}.$
\item $A_{3,5}=\begin{pmatrix}0&0&0&0&0&0&0&0&0\\ 0&0&0&0&0&0&0&0&0\\ 0&1&0&-1&0&1&0&0&0\end{pmatrix}.$
\item $A_{3,6}=\begin{pmatrix}0&0&0&0&0&0&0&0&0\\ 0&0&0&0&0&0&0&0&0\\ 0&1&0&1&0&1&0&0&0\end{pmatrix}.$
  \end{enumerate}
\end{itemize}
\emph{}\\

\item  Comparison with the list of \cite{GR}.
\emph{}\\
{\footnotesize{\begin{center}	
\begin{tabular}{|c|c|c|c|c|c|c|c|c}
 \hline
Algebra            & {Algebra in}               &Base change $g$& Isomorphism\\
from \cite{GR}&   this paper           &&  \\ \hline
$A_{3,2}$   & $As_{0}^{4}(3)(0)$                 &=& $A_{3,2}=As_{0}^{4}(3)(0)$\\
 $A_{3,3}^\delta$  &  $ As_0^4(3)(t), \ t=\delta \neq 0$&$ \begin{pmatrix}1&0&0\\ 0&0&1\\ 0&1&0\end{pmatrix}$& $A_{3,3}^t=g^{-1} As_0^4(3)(t)g^{\otimes 2}$  \\
 $A_{3,4}^\delta$   &     $As_{0}^{5}(3)(t), \ t=\delta$     &$ \begin{pmatrix}1&1&0\\ 0&0&1\\ 0&-1&0\end{pmatrix}$& $A_{3,4}^\delta=g^{-1}As_{0}^{5}(3)(t)(g^{\otimes 2}$  \\
 $A_{3,5}$  & $As_{0}^{1}(3)(-1)$           &$ \begin{pmatrix}0&0&1\\ 0&1&0\\ 1&0&0\end{pmatrix}$& $A_{3,5}=g^{-1}As_{0}^{1}(3)(-1)g^{\otimes 2}$ \\
$A_{3,6}$  &    $As_{0}^{1}(3)(1)$       &$\begin{pmatrix}0&0&1\\ 0&1&0\\ 1&0&0\end{pmatrix}$&$A_{3,6}=g^{-1}As_{0}^{1}(3)(1))g^{\otimes 2}$ \\
   \hline
\end{tabular}
\end{center}}}
\emph{}\\
\item \textbf{Conclusion}. The algebras $A_0^2(3)$ and
$A_0^3(3)$ are omitted.
\end{itemize}

\section{List of three-dimensional permutative algebras}
\begin{de} A permutative algebra $\mathbb{A}$ over a field
$\mathbb{F}$ is a vector space $\mathbb{V}$ with a $\mathbb{F}$-bilinear function $\cdot :\mathbb{V}\times \mathbb{V} \longrightarrow \mathbb{V}$ satisfying the identity
\begin{equation} \label{PA}
\begin{array}{rllllll}
(\mathrm x \cdot \mathrm y)\cdot \mathrm z &=& \mathrm x \cdot (\mathrm y \cdot \mathrm z)&=& \mathrm x \cdot (\mathrm z \cdot \mathrm y)
\ \text{(left perm-algebra)}\\
(\mathrm x \cdot \mathrm y)\cdot \mathrm z &=& (\mathrm y \cdot \mathrm x)\cdot \mathrm z&=& \mathrm x \cdot (\mathrm y \cdot \mathrm z) \ \text{(right perm-algebra)}.
\end{array}
\end{equation} \end{de}


This is the algebraic structure governed by the operad
$\mathbf{Perm}$ introduced in  \cite{ChL}. They are Koszul dual to the Pre-Lie operad.

Evidently, a commutative-associative algebra is permutative. A permutative algebra is associative, but an associative algebra
\begin{equation}\label{LPA}
\text{is left permutative if}\ \mathrm{x}(\mathrm{y}\mathrm{z}-\mathrm{z}\mathrm{y}) = 0,
\quad \mbox{for all }\ \mathrm{x},\mathrm{y},\mathrm{z}\in \mathbb{A},
\quad \mbox{equivalently,}\
[\mathbb{A},\mathbb{A}] \subseteq \mathrm{Ann}_l(\mathbb{A});
\end{equation}

\begin{equation}\label{RPA}
\text{is right permutative if}\ (\mathrm{x}\mathrm{y}-\mathrm{y}\mathrm{x})\mathrm{z} = 0,
\quad \mbox{for all }\ \mathrm{x},\mathrm{y},\mathrm{z}\in \mathbb{A},
\quad \mbox{equivalently,}\
[\mathbb{A},\mathbb{A}] \subseteq \mathrm{Ann}_r(\mathbb{A}).
\end{equation}
where $[\mathbb{A},\mathbb{A}]$, $\mathrm{Ann}_l(\mathbb{A})$, and $\mathrm{Ann}_r(\mathbb{A})$ are the commutant, left annihilator, and right annihilator, respectively.
%
%

B. Hou studied the extending structures problem of right-perm algebras and right-perm bialgebras from the
cohomological point of view (see \cite{Hou1}). The relationship
between the right permutative and other classes of algebras were given in \cite{Hou2}. Note that every metabelian Lie algebra can be embedded into
right permutative algebra under the commutator. The paper \cite{Ivan} contains the classification of three-dimensional complex left permutative algebras in connection with the algebraic and geometric classification problem of certain classes of algebras.

 Here are algebras in the list of Theorem \ref{3-dim} and the appearance of the algebras listed in Theorems 5 and 10 of \cite{Ivan}:
\begin{enumerate}
\item Commutative case (Theorem 5):

$\begin{array}{llllllll}
\text{Algebras in this paper}                                                          &\ \ \text{Algebras in}\ [1]             &\ \ \text{Algebras in}\ [11]\\
As_{1,1}^{1}(3)(t), \ t\in \mathbb{F} \ (\text{if}\ \mathbb{F}=\mathbb{C}\ \text{then}\ t \in \{0 ,1\})  &\ \begin{array}{llll} \mathbf{A}_{01}\cong A_{1,1}^1(3)(1)\\ \mathbf{A}_{02}= A_{1,1}^1(3)(0)\end{array}                                                                     &\begin{array}{llll} S_3^3 \cong A_{1,1}^1(3)(1)\\ S_4^3 = A_{1,1}^1(3)(0)\end{array} \\
As_{1,1}^{2}(3)(t), \ t\in \mathbb{F} \ (\text{if}\ \mathbb{F}=\mathbb{C}\ \text{then}\ t \in \{0 ,1\})  &\ \begin{array}{llll} \mathbf{A}_{07}\cong A_{1,1}^2(3)(1)\\ A_{1,1}^2(3)(0) \ \text{is omitted}\end{array}                                                                  &\begin{array}{llll} U_2^3\cong A_{1,1}^2(3)(1)\\ A_{1,1}^2(3)(0) \ \text{is omitted}\end{array}\\
As_{1,1}^{3}(3)                                                           &\quad \mathbf{A}_{10}\cong As_{1,1}^{3}(3)             &\ \ U_4^3\cong As_{1,1}^{3}(3)\\
As_{1,1}^{5}(3)                                                           &\quad \mathbf{A}_{09}= As_{1,1}^{5}(3)                 &\ \ U_0^3\cong As_{1,1}^{5}(3)\\
As_{1,1}^{6}(3)(-1)                                                       &  \quad omission                                                     &\quad omission\\
As_{1,1}^{7}(3)(-1)                                                       &   \quad omission                                                    &\quad omission\\
As_{1,1}^{8}(3)(t), \ t\in \mathbb{F} \ (\text{if}\ \mathbb{F}=\mathbb{C}\ \text{then}\ t \in \{0 ,1\})  &  \quad omission                      &\quad omission\\
As_{1,1}^{9}(3)(t), \ t\in \mathbb{F} \ (\text{if}\ \mathbb{F}=\mathbb{C}\ \text{then}\ t \in \{0 ,1\})  &   \quad omission                     &\quad omission\\
As_{1,1}^{10}(3)                                                          &\quad \mathbf{A}_{03}=As_{1,1}^{10}(3)                 &\ \ W_4^3 = As_{1,1}^{10}(3)\\
As_{1,1}^{11}(3)                                                          &\quad \mathbf{A}_{08}\cong As_{1,1}^{11}(3)            &\ \ U_3^3\cong As_{1,1}^{11}(3)\\
As_{1,1}^{13}(3)                                                          &\quad \mathbf{A}_{11}\cong As_{1,1}^{13}(3)            &\ \ S_2^3 \cong As_{1,1}^{13}(3)\\
As_{1,1}^{14}(3)                                                          &   \quad omission                                                    &\quad omission\\
As_{1,1}^{17}(3)(t), \ t\in \mathbb{F} \ (\text{if}\ \mathbb{F}=\mathbb{C}\ \text{then}\ t \in \{0 ,1\}) &   \quad omission                     &\quad omission\\
As_{0}^{1}(3)(1)                                                          &\quad \mathbf{A}_{05}\cong As_{0}^{1}(3)(1)            &\ \ W_3^3 \cong As_{0}^{1}(3)(1)\\
As_{0}^{3}(3)                                                             &\quad \mathbf{A}_{06}\cong As_{0}^{3}(3)               &\ \ S_1^3 \cong As_{0}^{3}(3)\\
As_{0}^{4}(3)(t), \ t\in \mathbb{F} \ (\text{if}\ \mathbb{F}=\mathbb{C}\ \text{then}\ t \in \{0 ,1\})    &\ \begin{array}{llll} \mathbf{A}_{04}= A_{0}^4(3)(0)\\ A_{0}^4(3)(1) \ \text{is omitted}\end{array} &\begin{array}{llll} W_1^3 \cong A_{0}^4(3)(0)\\ W_3^3(0) \cong A_{0}^4(3)(1) \end{array}\\
\end{array}$

\item Noncommutative case (Theorem 10):

$\begin{array}{llll}
\text{Algebras in this paper}                                                          &\quad \quad\quad \quad\text{Algebras in}\ [1]             &\ \text{Algebras in}\ [11]\\
As_2^2(3)(0)                                                                          &\quad \quad\quad \quad omission  &  \ U_1^3 \cong As_2^2(3)(0)
\\
As_2^4(3)                                                                              &\quad \quad\quad \quad omission   &\ omission 
\\
 As_2^5(3)                                                     &\quad \quad\quad \quad omission         & \ omission 
 \\
As_2^6(3)                                                      &\quad \quad\quad \quad \mathbf{A}_{23}\cong As_{2}^6(3)&\ W_7^3 = As_2^6(3) \\
 As_2^7(3)                                                     &\quad \quad\quad \quad \mathbf{A}_{24}\cong As_{2}^7(3)&\ W_8^3 =  As_2^7(3)\\
As_2^8(3)                                                      &\quad \quad\quad \quad\mathbf{A}_{20}\cong As_2^8(3)&\ omission \\
As_{1,3}(3)
		                                                       &\quad \quad\quad \quad\mathbf{A}_{17}= As_{1,3}(3)&\ C_3^3 \cong As_{1,3}(3)\\
As_{1,1/3}(3)
		                                                       &\quad \quad\quad \quad\mathbf{A}_{18}\cong As_{1,1/3}(3)&\ C_4^3 \cong As_{1,1/3}(3)\\
 As_{1,2}(3)                                                   &\quad \quad\quad \quad\mathbf{A}_{21}\cong As_{1,2}(3)&\ W_5^3 \cong  As_{1,2}(3) \\
As_{1,1/2}(3)
		                                                       &\quad \quad\quad \quad\mathbf{A}_{22}\cong As_{1,1/2}(3)&\ W_6^3 \cong As_{1,1/2}(3)\\
\end{array}$

$\begin{array}{llll}
   As_{1, 3/2}(3)
		                                                       &\quad \quad\quad \quad\mathbf{A}_{15}\cong As_{1,3/2}(3)&\ W_{10}^3 \cong  As_{1, 3/2}(3)\\
As_{1,2/3}(3)
		                                                       &\quad \quad\quad \quad\mathbf{A}_{16}\cong As_{1,2/3}(3)&\ W_9^3 \cong As_{1,2/3}(3)\\
        As_{1,1}^{12}(3)
		                                                       &\quad \quad\quad \quad\mathbf{A}_{19}\cong As_{1,1}^{12}(3)&\ C_2^3 \cong As_{1,1}^{12}(3)\\
As_0^{1}(3)(-1)                                                &\quad \quad\quad \quad\mathbf{A}_{12}\cong As_0^{1}(3)(-1)&\ C_1^3 \cong As_0^{1}(3)(-1)\\
As_0^{2}(3)                                                    &\quad \quad\quad \quad\mathbf{A}_{13} = As_0^{2}(3)&\ omission \\

As_0^{5}(3)(t), \ t\in \mathbb{F}
                              &\ \ \quad\quad \quad\begin{array}{llll} \textbf{A}^1_{14} \cong As_0^{5}(3)(0)\\ \mathbf{A}_{14}^\alpha\cong As_0^{5}(3)(\frac{\alpha^2-1}{4\alpha^2})\end{array}\ &\begin{array}{llll} W_2^3 \cong As_0^{5}(3)(0)\\ W_3^3(k) \cong As_0^{5}(3)(\frac{1}{k^2})\end{array}\\
\end{array}$
\end{enumerate}

Now from Theorem \ref{3-dim} verifying the conditions (\ref{LPA}) and (\ref{RPA}) we obtain the following classification result on three-dimensional permutative algebras over a field $\mathbb{F}$ of characteristic not two and three.
\begin{corollary} Every three-dimensional permutative algebra over a field $\mathbb{F}$ of characteristic not two and three is isomorphic to one the following pairwise non-isomorphic such algebras given by their tables of multiplications as follows
\footnotesize{\begin{center}	
\begin{tabular}{|c|c|l|c|c|c|c|c|c}
 \hline
N           & Notation
             &Table of multiplication& Type&$\begin{array}{llllll}\text{Algebra}\\\ \text{in}\ [1]\end{array} $  \\ \hline
$1$   & $As_2^6(3):$                 &$e_1e_1=e_1,\ e_3e_2=e_2,\ e_3e_3=e_3$& Right-perm&
\begin{tikzpicture}[baseline=(current bounding box.center)]
\draw (0,0) rectangle (0.5,0.5);
\draw (0,0) -- (0.5,0.5);
\draw (0,0.5) -- (0.5,0);
\end{tikzpicture}
\\  \hline
 $2$  &   $As_2^7(3):$&$ e_1e_1=e_1,\ e_2e_3=e_2,\ e_3e_3=e_3$& Left-perm& $\mathbf{A}_{24}$ \\ \hline
 $3$   &    $As_{1,3}(3)$    &$ e_1e_1=e_1,\ e_1e_2=e_2,\ e_1e_3=e_3$& Right-perm& \begin{tikzpicture}[baseline=(current bounding box.center)]
\draw (0,0) rectangle (0.5,0.5);
\draw (0,0) -- (0.5,0.5);
\draw (0,0.5) -- (0.5,0);
\end{tikzpicture} \\ \hline
 $4$   &     $As_{1,1/3}(3)$     &$ e_1e_1=e_1,\ e_2e_1=e_2,\ e_3e_1=e_3$& Left-perm& $\mathbf{A}_{18}$ \\ \hline
  $5$   &    $As_{1,2}(3)$    &$ e_1e_1=e_1, \ e_1e_3=e_3$& Both& $\mathbf{A}_{21}$ \\ \hline
  $6$   &     $As_{1,1/2}(3)$     &$ e_1e_1=e_1,\ e_3e_1=e_3$& Left-perm&  $\mathbf{A}_{22}$\\ \hline
$7$   &     $As_{1, 3/2}(3)$     &$ e_1e_1=e_1,\ e_1e_2=e_2,\ e_1e_3=e_3,\ e_2e_1=e_2$& Right-perm& \begin{tikzpicture}[baseline=(current bounding box.center)]
\draw (0,0) rectangle (0.5,0.5);
\draw (0,0) -- (0.5,0.5);
\draw (0,0.5) -- (0.5,0);
\end{tikzpicture} \\ \hline
 $8$  & $As_{1,2/3}(3)$           &$ e_1e_1=e_1,\ e_1e_2=e_2,\ e_2e_1=e_2,\ e_3e_2=e_3$& Left-perm& $\mathbf{A}_{16}$\\ \hline
 $9$  &    $As_{1,1}^1(t),\ t\in \mathbb{F}$       &$ e_1e_1=e_1,\ e_1e_2=e_2,\ e_2e_1=e_2,\ e_2e_2=te_1$&Both& $\begin{array}{ccc}\mathbf{A}_{01},\\ \mathbf{A}_{02}\end{array}$\\ \hline
  $10$  &    $As_{1,1}^2(t),\ t\in \mathbb{F}$       &$\begin{array}{ccc}  e_1e_1=e_1,\ e_1e_2=e_2,\\ e_2e_1=e_2,\ e_2e_2=te_1,\ e_3e_3=e_3\end{array} $&Both& $\mathbf{A}_{07}$\\ \hline
$11$  &    $As_{1,1}^3(3)$       &$\begin{array}{llllll} e_1e_1=e_1,\ e_1e_2=e_2,\ e_1e_3=e_3,\\ e_2e_1=e_2,\ e_3e_1=e_3,\ e_3e_3=e_3\end{array} $&Both&$\mathbf{A}_{10}$ \\ \hline
$12$  &    $As_{1,1}^5(3)$       &$ e_1e_1=e_1,\ e_1e_2=e_2,\ e_1e_3=e_3, \ e_2e_1=e_2,\ e_3e_1=e_3$&Both& $\mathbf{A}_{09}$ \\ \hline
$13$  &    $As_{1,1}^6(3)(-1)$       &$\begin{array}{llllll} e_1e_1=e_1,\ e_1e_2=e_2,\ e_1e_3=e_3, \ e_2e_1=e_2,\ e_2e_3=e_2,\\ e_3e_1=e_3,\ e_3e_2=e_2,\ e_3e_3=-e_1+e_2+2e_3 \end{array}$&Both& omitted\\ \hline
$14$  &    $As_{1,1}^7(3)(-1)$       &$\begin{array}{llllll} e_1e_1=e_1,\ e_1e_2=e_2,\ e_1e_3=e_3, \ e_2e_1=e_2,\ e_2e_3=e_2,\\ e_3e_1=e_3,\ e_3e_2=e_2,\ e_3e_3=-e_1+2e_3\end{array}$&Both &omitted\\ \hline
$15$  &    $As_{1,1}^8(3)(t),\ t\in \mathbb{F}$       &$\begin{array}{llllll} e_1e_1=e_1,\ e_1e_2=e_2,\ e_1e_3=e_3, \\ e_2e_1=e_2,\ e_2e_2=e_2,\ e_2e_3=e_3,\\ e_3e_1=e_3,\ e_3e_2=e_3,\ e_3e_3=te_1+te_2+e_3\end{array}$&Both&omitted \\
  &           &$\left(\begin{array}{llllll}As_{1,1}^8(3)(t)\cong As_{1,1}^8(3)(s)\\
  \begin{array}{llll}\text{if there exists}\ a\neq -\left(8t+\frac{3}{4}\right)\ \text{such that}\\ a^2=1+64st+8(s+t),\ \text{where}\ s\neq-\frac{1}{8}\end{array}\end{array}\right)$&&\\ \hline
  $16$  &    $As_{1,1}^9(3)(t),\ t\in \mathbb{F}$       &$\begin{array}{llllll} e_1e_1=e_1,\ e_1e_2=e_2,\ e_1e_3=e_3, \\ e_2e_1=e_2,\ e_2e_2=e_2,\ e_2e_3=e_3,\\ e_3e_1=e_3,\ e_3e_2=e_3,\ e_3e_3=te_1+te_2\end{array}$ &Both&omitted\\
&           & $As_{1,1}^9(3)(t)\cong As_{1,1}^9(3)(a^2t)$, where $a \in \mathbb{F}^*$ & &\\ \hline
  \end{tabular}
\end{center}}

\small{\begin{center}	
\begin{tabular}{|c|c|l|c|c|c|c|c|c}
 \hline
$17$  &    $As_{1,1}^{10}(3)$       & $e_1e_1=e_1$&Both& $\mathbf{A}_{03}$\\ \hline
$18$  &    $As_{1,1}^{11}(3)$       & $e_1e_1=e_1, \ e_2e_3=e_2,\ e_3e_2=e_2,\ e_3e_3=e_3$&Both& $\mathbf{A}_{08}$\\ \hline
$19$  &    $As_{1,1}^{13}(3)$       &$e_1e_1=e_2,\ e_3e_3=e_3$ &Both&$\mathbf{A}_{11}$ \\ \hline
$20$  &    $As_{1,1}^{14}(3)(t),\ t\in \mathbb{F}$       &$e_2e_3=e_2,\ e_3e_2=e_2, \ e_3e_3=e_1+te_2+e_3$ &Both&omitted \\ \hline
$21$  &    $As_{1,1}^{15}(3)$       & $e_2e_3=e_2,\ e_3e_2=e_2, \ e_3e_3=e_2+e_3$&Both&omitted \\ \hline
$22$  &    $As_{0}^{1}(3)(-1)$       & $e_2e_3=-e_1,\ e_3e_2=e_1$&Right-perm&\begin{tikzpicture}[baseline=(current bounding box.center)]
\draw (0,0) rectangle (0.5,0.5);
\draw (0,0) -- (0.5,0.5);
\draw (0,0.5) -- (0.5,0);
\end{tikzpicture}\\ \hline
$23$  &    $As_{0}^{2}(3)$       & $e_1e_1=e_2,\ e_1e_3=e_2,\ e_3e_1=-e_2$&Both &$\mathbf{A}_{13}$\\
   \hline
$24$  &    $As_{0}^{3}(3)$       &$e_2e_3=e_1,\ e_3e_2=e_1,\ e_3e_3=e_2$ &Both&$\mathbf{A}_{12}$ \\ \hline
$25$  &    $As_{0}^{4}(3)(t),\ t\in \mathbb{F}$       &$\begin{array}{llllll}e_1e_1=e_2,\ e_3e_3=te_2 \\ (As_{0}^{4}(3)(t)\cong As_{0}^{4}(3)(a^2t), \ a\in \mathbb{F}^*)\end{array}$ &Both &$\mathbf{A}_{06}$\\ \hline
$26$  &    $As_{0}^{5}(3),\ t\in \mathbb{F}$       & $\begin{array}{llllll}e_1e_1=e_2,\ e_3e_1=e_2,\ e_3e_3=te_2\\
(As_0^{5}(3)(t)\simeq As_0^{5}(3)(-t),\\ \text{if}\ u^2=-1\ \text{has a solution in}\ \mathbb{F})\end{array}$&Both &$\mathbf{A}_{14}$\\
   \hline
\end{tabular}
\end{center}}
\emph{}\\

In the table $\begin{tikzpicture}[baseline=(current bounding box.center)]
\draw (0,0) rectangle (0.5,0.5);
\draw (0,0) -- (0.5,0.5);
\draw (0,0.5) -- (0.5,0);
\end{tikzpicture}$ means not applicable.
\end{corollary}



\begin{thebibliography}{9}
\bibitem{Ivan} {\sc H. Abdelwahab, I. Kaygorodov, R. Lubkov}, The algebraic and geometric classification of
right alternative and semi-alternative algebras, 2025, https://doi.org/10.48550/arXiv.2505.00720v2.
\bibitem{B} {\sc U. Bekbaev}, Classification of two-dimensional algebras over any base field, 2023, {\it AIP Conference Proceedings}, 2880, 030001, https://doi.org/10.1063/5.0165726.
\bibitem{ChL} {\sc  F. Chapoton, M. Livernet,}  Pre-Lie algebras and the rooted trees operad, 2001, {\it International Mathematics Research Notices}, 8, 395--408.
 \bibitem{FP} {\sc A. Fialowski, M. Penkava,}   The Moduli Space of 3-Dimensional Associative Algebras, Communications in Algebra, 2009, 37(10), 3666--3685.
 \bibitem{FPP} {\sc A. Fialowski, M. Penkava, M. Phillipson},   Deformations Of Complex 3-Dimensional Associative Algebras, Journal of Generalized Lie Theory and Applications, 2011, 1--27, DOI: 10.4303/jglta/G110102.
\bibitem{Gabriel}{\sc P. Gabriel}, Finite representation type is open, \emph{Proceedings of}
ICRAI, \emph{Ottawa} 1974, \emph{Lect. Notes in Math.}, 1975, 488, Springer-Verlag.
\bibitem {GR}  {\sc W.A. De Graaf,} Classifification of nilpotent associative algebras of small dimension, International Journal of Algebra and Computation, 2018, 28(1), 133--161.
\bibitem {GM} {\sc M. Goze, A. Makhlouf}, Classification of rigid associative algebras in low dimensions, in: Lois d’Alg\'{e}bres et Variét\'{e}s Alg\'{e}briques (Colmar, 1991), in: Travaux en Cours, vol. 50, Hermann, Paris, 1996, pp. 5--22.
\bibitem {Hou1} {\sc B. Hou}, Extending structures for perm algebras and perm
bialgebras, 2024, {\it Journal of Algebra}, 2024, 649, 392--432. https://doi.org/10.1016/j.jalgebra.2024.03.018
\bibitem {Hou2} {\sc B. Hou}, Affinization of Dendriform $D$-bialgebras, Lie bialgebras and
solutions of classical Yang-Baxter equation, 2026, https://doi.org/10.48550/arXiv.2601.17456
\bibitem {KSTT} {\sc Y. Kobayashi, K. Shirayanagi, M. Tsukada, S. Takahasi}, A complete classification of three-dimensional algebras
over $\mathbb{R}$ and $\mathbb{C}$, Asian-European Journal of Mathematics, 2021, 14(8) 2150131, DOI: 10.1142/S179355712150131X.
\bibitem{M}{\sc G. Mazzola}, The algebraic and geometric classification of
associative algebras of dimension five, \emph{Manuscripta Math.}, 1979, 27, 1--21.
\bibitem{P}{\sc B. Poonen}, Isomorphism types of commutative algebras of finite rank over an algebraically closed field, 2008, In Computational Arithmetic Geometry, Contemporary Mathematics (American Mathematical Society Providence),  111--120.
\bibitem {R} {\sc I.S.Rakhimov}, Algebraic Structures on Two-Dimensional Vector Space Over Any Basic Field, {\it Malaysian Journal of Mathematical Sciences}, 18(2) (2024), 227--257, doi:10.47836/mjms.18.2.02.
\bibitem{S}{\sc B.G. Scorza}, Le algebre del 3 ordine, \emph{Atti. Nap.}, 1938, 20(13), 20(14) (in Italian).
\end{thebibliography}
 	\end{document}